\pdfoutput=1
\documentclass[a4paper,legno,11pt]{article}

\usepackage[accsupp]{axessibility}

\usepackage[utf8]{inputenc}
\usepackage[T1]{fontenc}
\usepackage{xspace,
	amsthm,
	amsmath,
	amssymb,
	bbm,
	tikz,
	graphicx,
	subfig,
	keyval,
       xcolor,
  pgfplots}

\DeclareGraphicsExtensions{.png, .pdf}

\usepackage[hidelinks]{hyperref}

\usepackage{geometry}
\tikzstyle{nodino}=[circle,draw,fill,inner sep=0pt,minimum size=0.5mm]
\tikzstyle{infinito}=[circle,inner sep=0pt,minimum size=0mm]
\tikzstyle{nodo}=[circle,draw,fill,inner sep=0pt, minimum size=0.5*width("k")]
\tikzstyle{nodo_vuoto}=[circle,draw,inner sep=0pt, minimum size=0.5*width("k")]
\usetikzlibrary{graphs}
\tikzset{every loop/.style={min distance=10mm,in=300,out=240,looseness=10}}
\tikzset{place/.style={circle,thick,draw=blue!75,fill=blue!20,minimum
		size=6mm}}
\tikzset{place2/.style={circle,thick,draw=red!75,fill=red!20,minimum
		size=6mm}}

\newcommand{\R}{{\mathbb R}}

\newcommand{\G}{{\mathcal{G}}}

\newcommand{\E}{{\mathcal{E}}}
\newcommand{\F}{{\mathcal{F}}}

\newcommand{\Htau}{H_\tau^1(\mathbb{R})}

\newcommand{\uLtwor}{\|u\|_{L^2(\mathbb{R})}}

\newcommand{\V}{\mathcal{V}}
\newcommand{\om}{\omega}

\theoremstyle{plain} 
\newtheorem{thm}{Theorem}[section]

\newtheorem{prop}[thm]{Proposition} 

\theoremstyle{definition}

\theoremstyle{definition}

\theoremstyle{remark}

\title{Non-Kirchhoff Vertices and Nonlinear Schr\"{o}dinger Ground States on graphs\footnote{We acknowledge that the present research was partially supported by MIUR Grant - Dipartimenti di Eccellenza 2018-2022 n. E11G18000350001.}}

\author{Riccardo Adami$^\dagger$, Filippo Boni$^{\dagger,\star}$, Alice Ruighi$^{\dagger,\star}$
\\ \ \\{\small  {$^\dagger$}Dipartimento di Scienze
Matematiche ``G.L. Lagrange'', Politecnico di Torino } \\ {\small
Corso Duca degli Abruzzi, 24, 10129 Torino, Italy} \\ \ \\
{\small  {$^\star$}Dipartimento di Matematica ``G. Peano'',  Universit\`a di Torino } \\ {\small
Via Carlo Alberto 10, 10123 Torino, Italy
}}

\date{}
\begin{document}
\maketitle

\begin{abstract}
We review some recent results and announce some new ones on the problem of the existence of ground states for the Nonlinear Schr\"odinger Equation on graphs endowed with vertices where the matching condition, instead of being free (or Kirchhoff's), is non-trivially interacting. In this category fall Dirac's delta conditions, delta prime, F\"ul\"op-Tsutsui, and others.
\end{abstract}

\section{Introduction}
In this paper we aim at presenting some new trends in the study of the Nonlinear Schr\"odinger equation on metric graphs. Specifically, we focus on cases where non-Kirchhoff's conditions are imposed at the vertices. 

Roughly speaking, metric graphs, also called networks,
 are a simple mathematical scheme used to build up models of branched structures, as they are made of one-dimensional segments, named {\em edges} or {\em arcs}, that can intersect in points called {\em vertices}. The only geometric information available for such structures is the length of the edges, whereas no curvature, angles, nor embedding in higher dimensional settings are considered. Therefore, metric graphs are suitable to provide simple models for the propagation of signals of various origin: acoustic, electric and also matter waves (for an exhaustive physical introduction see for instance \cite{noja14}).

As it is intuitive, in order to be defined, a wave dynamics on a network requires an evolution equation, describing the change of the profile in time inside the edges, together with a matching condition at the vertices, that rules the transmission and the reflection of the signals 
  when crossing a junction. In the context of {\em quantum graphs}, namely metric graphs in which the ruling equation is the linear Schr\"odinger, the problem of finding all admissible transmission and reflection rates, maps to the issue of finding all possible self-adjoint extensions of the restriction of the Laplacian to functions that vanish in a neighbourhood of every vertex. This problem was exhaustively treated in a seminal paper by Kostrykin and Schrader \cite{kost}. In fact, self-adjointness is the translation in the language of operators of the consevation of the total probability, that is a fundamental and inescapable requirement for every genuine quantum theory. In the context of quantum graphs, self-adjointness, and therefore conservation of probability, is encoded exactly in the law ruling the transmission and the reflection at the vertices. More generally, the task of finding all self-adjoint extensions of a symmetric operator is the key ingredient of the definition of point interactions, i.e. potentials located at a single point in space. For a collection of fundamental results we refer the reader to the treatise \cite{alb}. On the other hand, for the application of the theory of self-adjoint extensions to graphs, the basic reference is the monograph \cite{kuchment}. Among all the possible matching conditions, the most used are those named after Kirchhoff as a reminescence of Kirchhoff's law for linear circuits. Such conditions prescribe that the sum of the derivatives of the wave function ingoing into every vertex equals zero (for the mathematical formulation of such conditions see \eqref{kirch}). In the particular case of a vertex attached to one edge only, Kirchhoff's condition reduces to Neumann's, while in the case of two edges only concurring to the same vertex, Kirchhoff's condition restores the requirements of continuity and differentiability at the point occupied by the vertex. Furthermore, Kirchhoff's conditions naturally arise when dealing with the search for ground states, namely with the minimization of the energy functional under, as a unique constraint, the value of the mass. Finally, when imposing Kirchhoff's conditions, transmission and reflection coefficients are independent of the momentum.
For all these reasons, Kirchhoff's conditions have been assumed as the most natural, and hence the most widely studied. However, it is not clear whether in relevant physical contexts such conditions prove able to fit the dynamics well. On the contrary, F\"ul\"op, Tsutsui and Cheon (\cite{ft,ftc}) suggested that some other conditions could be more satisfactory from the point of view of invariance laws. Furthermore, the presence of non-trivial, localized interactions near junctions make quite unlikely that Kirchhoff's conditions would supply a predictive model. Therefore, quite recently a programme for the study of nonlinear dynamics on graphs including non-Kirchhoff's conditions has started. Here we report on the first results and perspectives of this project with a particular emphasis on the existence and stability of ground states for the Nonlinear Schr\"odinger equation on graphs, endowed with non-Kirchhoff's conditions at the vertices.

In the literature, metric graphs are graphs $\G=(\mathcal{V},\E)$, where $\mathcal{V}$ is the set of the vertices and each edge $e \in \E$ is identified either with a bounded and closed interval $I_e=[0, l_e]$ or with a halfline $I_e=[0,+\infty)$.
\begin{center}\vspace{1cm}
\includegraphics[width=0.5\linewidth]{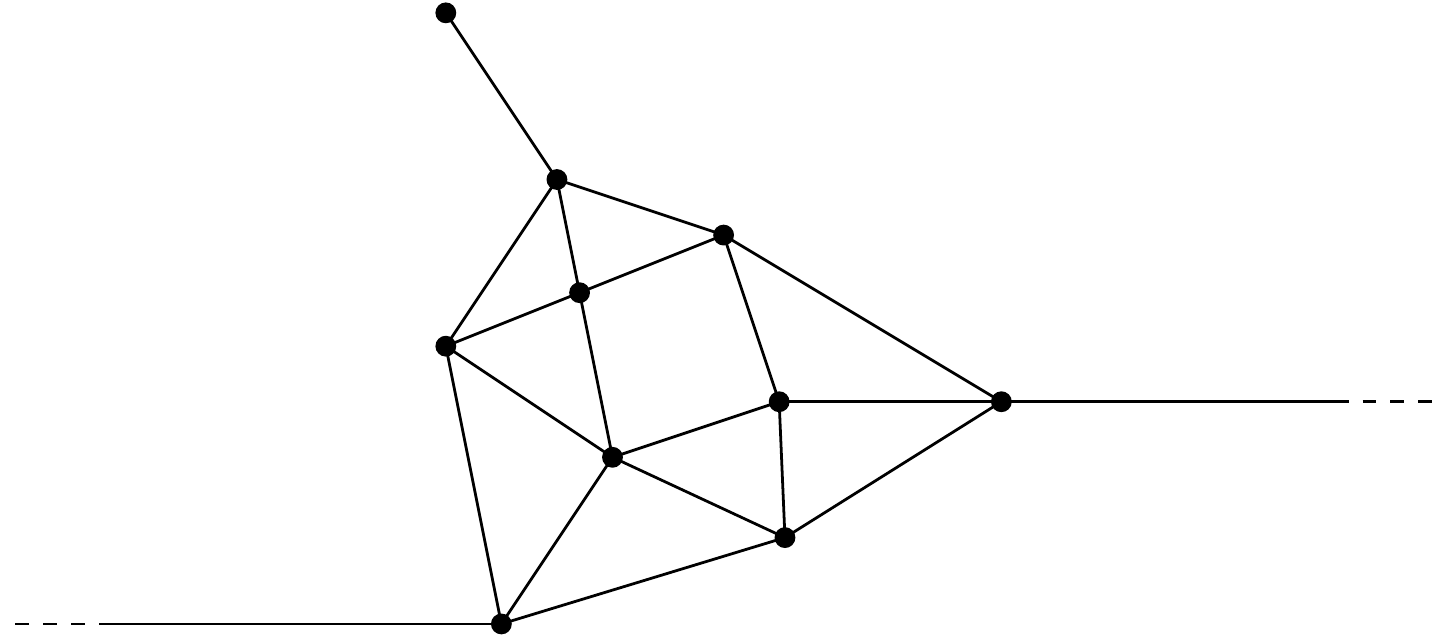}
\captionof{figure}{$\G=(\mathcal{V},\E)$}
\end{center}\vspace{1cm}
Owing to the metric structure, on this kind of graphs one can define real or complex-valued functions and the related ordinary functional spaces, like Lebesgue and Sobolev spaces, as $$L^p(\G)= \bigoplus_{e \in \mathcal{E}} L^p(I_e) $$ and $$H^1(\G)= \bigg\{ u \in \bigoplus_{e \in \mathcal{E}} H^1(I_e) \, : \,  u  \text{ is continuous at } v,  \, \forall v \in \mathcal{V}  \bigg\}.$$
Consequently, a function $\psi$ on a metric graph can be seen as a collection of functions $\psi=(\psi_e)_{e\in \mathcal{E}}$ and the variable $x$ on the graph runs through the collection of all variables $x_e$ defined on each edge.
\\
As already specified, dynamics on metric graphs can provide a good approximation for the evolution of systems located on ramified structures, i.e. systems for which locally there is a privileged direction for the propagation of signals, since the dimensions transverse to that of propagation are negligible compared to the longitudinal one. Such structures are often referred to as {\em quasi one-dimensional}.

 This line of investigation has recently undergone a dramatic boost due to its closeness to the problem of the evolution of Bose-Einstein Condensates (BEC).\\ Let us just briefly recall that a Bose-Einstein Condensate is a system composed by a large number of identical bosons (often alkali atoms) whose spatial confinement is usually realized by magneto-optical traps. Bose in 1924 \cite{bose} and Einstein in 1925 \cite{einstein} predicted that, under a critical value of the temperature, the state of the whole system collapses into a non-classical state in which each particle acquires the same wave function, called {\em wave function of the condensate}, that solves the following variational problem
\begin{equation}
    \label{GPvar}
    \min_{\substack{u \in H^1(\Omega), \\ \int |u|^2= N}} E_{GP}(u),
\end{equation}
where $\Omega$ is the trap in which the particles are confined, $N$ is the number of the particles of the system and finally $E_{GP}$ is the Gross-Pitaevskii functional defined as
\begin{equation*}
    \label{GP}
    E_{GP}(u)= ||\nabla u||^2_{L^2(\Omega)}+8 \pi \alpha ||u||^4_{L^4(\Omega)}.
\end{equation*}
In equation \eqref{GP}, $\alpha$ is the scattering length of the two-body interaction between the particles in the condensate. \\Provided that it exists, a solution $u_\om$ to the variational problem \eqref{GPvar} 
must satisfy the Euler-Lagrange equation
$$ - \Delta u+32 \pi \alpha |u|^{2} u + \omega u = 0,$$
where $\omega$ arises as a Lagrange multiplier and depends on $N$, and it is immediately seen that the function $\psi (t,x) = e^{i \om t} u_\om (x)$
is a solution, in particular a {\em standing wave}, to the Gross-Pitaevskii equation
\begin{equation*}
    i \partial_t \psi = -\Delta \psi +32 \pi \alpha |\psi|^{2} \psi.
\end{equation*}

The interaction between atoms in a BEC is usually repulsive, so that the sign of $\alpha$ is in general positive, but it is nowadays possible to tune such interaction through a mechanism called Feshbach resonance \cite{ccdt19}, so that it become possible to create collapsing condensates, by making $\alpha$ to become negative. This fact makes interesting to study Gross-Pitaevskii equation with a focusing nonlinearity. 
Concerning the trap $\Omega$, it is generically intended to be a smooth region of the three-dimensional space. It is then possible to arrange experimentally elongated traps, and also branched traps by the insertion of suitable junctions. In this cases, the trap remains genuinely three-dimensional, but
it is commonly accepted, even though a general and rigorous proof is still lacking, that, as a result of a suitable shrinking limit, Gross-Pitaevskii equation for a three dimensional system can be approximated by a  focusing nonlinear Schr\"{o}dinger equation on a metric graph
like
\begin{equation}
\label{partial}
i \partial_t \psi = H\psi - |\psi|^{p-2}\psi \quad \text{with} \quad p>2,
\end{equation}
where the graph can be understood as a quasi one-dimensional skeleton of the original elongated and branched trap.

In \eqref{partial}, $H$ is a self-adjoint extension of the Laplace operator 
\begin{equation}
\label{laplacian}
-\Delta: \bigoplus_{e \in \E} C^\infty_0\left(\mathring{I_e}\right)\to L^2(\G),   
\end{equation}
where $\mathring {I_e}$ is the {\em open} interval  corresponding to the edge $e$, and the space $C^\infty_0\left(\mathring{I_e}\right)$ collects all smooth and compactly supported functions on $\mathring{I_e}$.

A prolific research line is focused on the search for standing waves of \eqref{partial}, i.e. solutions of \eqref{partial} of the form $\psi(t,x)=e^{i \om t}\phi_{\om}(x)$, where $\omega\in\R$ and $\phi_\om$ solves the stationary equation
\begin{equation}
\label{stationary}
H\phi - |\phi|^{p-2}\phi+\om \phi=0.
\end{equation}
A significant part of the literature has dealt with stationary solutions for the focusing nonlinear Schr\"{o}dinger equation when Kirchhoff's boundary conditions are imposed at the vertices, namely
\begin{equation}
\label{stateq}
-\phi''- |\phi|^{p-2}\phi+\om \phi=0\quad\text{in}\quad \mathring{I_e},\quad \forall e\in \mathcal{E},
\end{equation}
coupled with the boundary conditions
\begin{equation} \label{kirch}
\begin{cases}
\phi_{e_1}(v)=\phi_{e_2}(v), \quad \forall e_1,e_2\succ v,\quad \forall v\in \V,\\
\sum_{e\succ v} \frac{d\phi_e}{dx_e}(v)=0,\quad \forall v\in \V,
\end{cases}    
\end{equation}
where $e\succ v$ means that $e$ is an incoming edge of the vertex $v$ and $\frac{d\phi_e}{dx_e}(v)$ stands for $\phi'_e(l_e)$ if $x_e=l_e$ at $v$ and for $\phi'_e(0)$ if $x_e=0$ at $v$.\\
In this context, two main variational approaches have been used to find solutions to this equation, both consisting in minimizing a proper functional under some additional constraints.
In the first approach one looks for critical points of the energy functional 
\begin{equation}
\label{energy}
 E(u,\G)= \frac{1}{2}||u'||^2_{L^2(\G)}-\frac{1}{p}||u||^p_{L^p(\G)}    
\end{equation}
on the manifold  
\begin{equation*}
    \label{mass}
    H^1_{\mu}(\G)=\left\{u\in H^1(\G):\int_{\G} |u|^2 dx = \mu\right\},
\end{equation*}
where $p >2$ and $\mu$ is a fixed parameter usually called the mass. In the second case, instead, one is interested in critical points of the action functional
\begin{equation}
    \label{action}
    S_{\om}(u,\G)=E(u,\G)+\frac{\om}{2}||u||^2_{L^2(\G)}
\end{equation}
under a constraint called Nehari's constraint, namely
\begin{equation*}
    I_{\om}(u,\G)=0,
\end{equation*}
where
\begin{equation*}
    \label{nehari}
    I_{\om}(u,\G)=||u'||^2_{L^2(\G)}-||u||^p_{L^p(\G)}+\om||u||^2_{L^2(\G)}.
\end{equation*}

\noindent
This constraint is considered as natural since it hosts all stationary points of the action functional \eqref{action}.

\noindent In the following, we will refer to global minimizers of \eqref{energy} or \eqref{action} as ground states, regardless of the functional they minimize.   
The difference between the two approaches reflects on the parameter $\om \in \R$ in the equation \eqref{stationary}. In fact, it can be an unknown of the problem and be interpreted as a Lagrange multiplier like in the former approach or it can be given, like in the latter. \\\\

The forerunner of the study of nonlinear evolution on metric graphs is considered to be a paper by Ali Mehmeti \cite{alimehmeti} that dates back to 1984, but it is in the last decade that the study of the NLSE with Kirchhoff's conditions has developed. Several papers \cite{ast15,ast-jfa,ast17} study the existence of ground states for the energy functional under the mass constraint. They analyse the problem for metric graphs with a finite number of vertices and at least one halfline, distinguishing between cases $p \in (2,6)$ (known as {\em subcritical}) and  $p=6$ (called the {\em critical} case). In the first case, the authors show how the topology of the graph can affect the existence of ground states, or, on the other hand, when it depends on the interplay between the metric features of the graph and the mass $\mu$. In the critical case, instead, the mass assumes a crucial role for the existence of ground states. In particular, differently from the case of the real line, the authors prove that ground states can exist not only for a critical value of the mass, but for a whole interval of masses. 
In addition, the existence of ground states for periodic graphs has been studied in \cite{adst,adr,dovetta-per,ad}, while the problem on infinite metric trees has been approached and partially solved in \cite{dst}. Finally, existence of ground and/or bound states for the NLS on graphs with a nonlinearity concentrated on a subgraph has been variously explored, for instance in \cite{serratentarellly,dt,t}. \\

On the other hand, the search for ground states for the NLS on the line with delta or delta-prime interactions, that can be thought of as a graph with a vertex and two infinite edges with a non-Kirchhoff's condition, is older than the extended specific research on graphs \cite{foo,an}. Point interactions on the real line have been previously studied in the time-dependent setting by Caudrelier et al. using integrability tools \cite{cmr1,cmr2}. We stress however that these results are bounded to the case of the cubic nonlinearity, i.e. $p=4$, that corresponds to the integrable case \cite{zakharov}. Other non-free boundary conditions on graphs have been extensively investigated in connection with the integrability features. A breakthrough result, due to Matrasulov and coworkers, is the discovery of a class of non-reflecting matching conditions 
that make the cubic NLS on graphs inherit the integrability from the corresponding one-dimensional system (\cite{ssms,sbmnu}) and, at least for star graphs, make possible to restore the structure and the methods typical of integrable systems, like Lax pairs and inverse scattering: this was accomplished in \cite{c15}, whose results extend to non-integrable boundary conditions too. Other milestone results on graphs with the same non-reflecting (hence non-Kirchhoff) conditions have been obtained by Pelinovsky and collaborators \cite{kp2,kp,kpg} in a series of works where the spectral stability of special solutions (like half-solitons or shifted states) was investigated.

\noindent The main purpose of this review is to present some of the results in the literature or still in preparation which involve nonlinear Schr\"odinger equations on metric graphs when non-Kirchhoff's conditions are imposed at the vertices.\\
The physical motivations for the introduction of such conditions at the vertices rely on the necessity to represent some inhomogeneity or defect in the medium in which the dynamics occurs. Rigorous studies of NLSE in presence of impurities described by point interactions have been given along several lines, with a special consideration for the so-called {\em{delta interaction}}.\\
As already specified, delta interactions have been the oldest non-Kirchhoff's conditions to be studied and mathematically they are described as conditions localized at the vertices $v \in \mathcal{V}$, involving both the value of the function and its derivative. Specifically, they are defined by 
\begin{equation}
    \label{deltaint}
      \begin{cases}
      \phi_{e_1}(v)=\phi_{e_2}(v), \quad \forall e_1, e_2\succ v,\quad \forall v\in \V \\
      \sum_{e\succ v} \frac{d\phi_e}{dx_e}(v) =\alpha \phi(v), \quad \forall v \in \V
\end{cases}    
\end{equation}
with $\alpha \in \R$, and they are obtained as the result of a proper self-adjoint extension of the Laplace operator in \eqref{laplacian}.\\ 
In constrast to Kirchhoff's conditions, up to now the study of delta conditions has been confined on simple graphs with a single vertex, namely the real line $\R$ or more generally star graphs $S_N$ with $N$ halflines.
 This is the starting point for possible future studies on more general graphs, a step that turns out to be highly non-trivial, considering that for general power nonlinearities it is not  possible to directly make recourse to abstract methods used e.g. in \cite{cr}.

\noindent Through the study of delta interactions one is naturally led to consider two different classes of non-Kirchhoff's conditions: the first one is the class of linear non-Kirchhoff's conditions, which ensure self-adjointness of the operator defined on the metric graph and include \emph{delta prime}, \emph{dipole} and \emph{F\"{u}l\"{o}p-Tsutsui's conditions}, while the second one is the class of \emph{nonlinear delta interactions}, obtained replacing the real number $\alpha$ with a nonlinear function of $\phi$. In the latter case, two nonlinearities coexist: the standard one, given by the $p$-th power of the $L^p$ norm, and a pointwise one.\\
The review is organized as follows: in Section 2 we will introduce the linear non-Kirchhoff's conditions, collecting some old results and presenting other ones more recent. In Section 3 some new results concerning the nonlinear delta conditions are shown.  

\section{Linear non-Kirchhoff's condition}

As outlined in the Introduction, linear non-Kirchhoff's conditions are a family of conditions imposed at the vertices of a metric graph $\G$ in such a way that the operator $H$ in \eqref{partial} turns out to be self-adjoint. Among the first ones who studied non-Kirchhoff's conditions, in the integrable cubic case there were Caudrelier et al. \cite{cmr2}, who presented a family of point interactions that preserves the quantum integrability and Goodman et al. \cite{ghw} and Holmer et al. \cite{hmz}, who introduced delta interactions and started the study of the existence and the stability of solutions of the NLSE.\\
As for the problem with Kirchhoff's conditions, two main approaches have been carried on: the first one concerns the minimization of the energy under the mass constraint, while the second one is based on the minimization of the action functional restricted to the Nehari's manifold.\\

\subsection{Minimization of the energy under the mass constraint}

This first approach has been used in \cite{anv}, where the authors study existence of 1D ground states and their orbital stability when a point interaction is placed at the origin of the real line and the standard nonlinearity is subcritical, i.e. when $2<p<6$. In particular, they are interested in three different conditions at the origin:
\begin{itemize}
    \item {attractive delta conditions, i.e.
    \begin{equation}
      \label{delta}
      \begin{cases}
      \phi(0_+)=\phi(0_-),  \\
      \phi'(0_-)-\phi'(0_+)=\alpha \phi(0),
      \end{cases}
      \end{equation}
      where $\alpha>0$,}
     \item{delta prime conditions, that are
      \begin{equation}
      \label{deltapr}
      \begin{cases}
      \phi'(0_+)=\phi'(0_-),\\
      \phi(0_-)-\phi(0_+)= \beta \phi'(0),
\end{cases}    
\end{equation}
where $\beta >0$,
      }
      \item{dipole conditions, i.e.
     \begin{equation}
      \label{dipole}
      \begin{cases}
      \phi(0_+)=\tau \phi(0_-),\\
      \phi'(0_-)- \tau \phi'(0_+)=0,
\end{cases}    
\end{equation}
with $\tau \in \R$.}
    \end{itemize}
We recall that these boundary conditions are induced by three of the possible self-adjoint extensions of the 1D Laplacian on the real line and up to now they are the most studied in the stationary setting. In this context, the authors of \cite{anv} prove an abstract theorem that revisits the concentration-compactness method by Cazenave and Lions \cite{lions} and which is suitable to treat all these three inhomogeneities. Applying this general result to all these three cases, it is straightforward that ground states exist and they are orbitally stable: moreover, thanks to the one dimensional structure, it is possible to compute explicitly solutions.\\
These results do not reveal any substantial novelty moving from the Kirchhoff's case to these non-Kirchhoff's conditions: indeed, a ground state exists for every mass also in the Kirchhoff's case.\\

\noindent The relevance of point interactions becomes clearer  passing from the real line $\R$ to star graphs $S_N$. Indeed, in the paper \cite{acfn2013}, the concentration-compactness method was adapted to the case of an attractive delta interaction localized at the vertex of a star graph and prove that there exists a threshold value of the mass $\mu^*$ such that, under this value, the ground state exists, is symmetric, decreasing on each halfline and orbitally stable: among all the stationary states, this is called the $N$ tail state (see Figure \ref{3star}). In particular, the result is valid both in the subcritical and in the critical case.\\
It is important to notice that, in the Kirchhoff's case, ground states do not exist for any value of the mass and for this reason the delta interaction is crucial to gain the existence of ground states for small masses.    Moreover, the $N$ tail state turns out to be orbitally stable for any value of the mass $\mu>0$ (see \cite{acfn2016}) : this means that the orbital stability stands even when the $N$ tail state is not a ground state but only a local minimizer of the energy functional.  

\begin{center}\vspace{1cm}
\begin{tikzpicture}[xscale= 0.7,yscale=0.5]
\node at (-9,0) [infinito] (GS) {};
\node at (6,3.5) [infinito] (GS1) {};
\node at (7,-2) [infinito] (GS2) {};
\node at (-1,0) [nodo] (G1) {};

\node at (-2,-2)  [infinito] (G) {$S_3$};

\draw[-] (GS)--(G1);
\draw[-] (G1)--(GS1);
\draw[-] (G1)--(GS2);
\draw[dashed] (GS)--(-11,0);
\draw[dashed] (GS1)--(8,4.5);
\draw[dashed] (GS2)--(9,-2.66);

\node[-] at (-1,4) [infinito] (F1) {};
\node[-] at (6,3.7) [infinito] (F2) {};
\node[-] at (7,-1.8) [infinito] (F3) {};

\draw[dashed] (G1)--(F1);
\draw[-,thick] (-8.55,.2) to [out= 5, in = 240] (F1);
\draw[-,thick] (F1)  to [out= -25, in = 200]  (F2);
\draw[-,thick] (F1)  to [out= -50, in = 170] (F3);
\end{tikzpicture}
\captionof{figure}{3-tail state}
\label{3star}
\end{center}\vspace{1cm}

\subsection{Minimization of the action under the Nehari's constraint}

As noticed in the introduction, there is a second way to approach the study of the NLS equation and it consists in minimizing a constrained action funcional. \\
  The usual aims in this setting are to identify stationary states, characterize ground states and show if standing waves are stable or not for any value of the power $p$. We remark that in this section stability always stands for orbital stability and the outcomes concerning it are achieved taking advantage of the well known theory by Grillakis, Shatah and Strauss \cite{gss1,gss2}.  \\\\
Fukuizumi et al. (see \cite{foo}, \cite{fj} and \cite{coz}) investigate the previous problems on the real line with a delta type defect located at the origin, analysing both the attractive case ($\alpha >0$) and the repulsive one ($\alpha <0$). In particular they use the conditions defined in \eqref{delta}
with $\alpha\in \R$ and study the existence and the stability of global minimizers of the action functional 
\begin{equation} 
\label{Somega}
S_{\om,\alpha}(u)=\frac{1}{2}||u'||^2_{L^2(\R)}-\frac{1}{p}||u||^p_{L^p(\R)}-\frac{\alpha}{2}|u(0)|^2+\frac{\om}{2}||u||^2_{L^2(\R)}
\end{equation}
under the Nehari's constraint. A first difference between the results obtained in the attractive case and in the repulsive one concerns the functional space in which the minimization holds. In fact, if for $\alpha >0$ the existence of the ground state is proved in $H^1(\R,\mathbb{C})$, for $\alpha < 0$ the same result is shown only on the subspace of the even functions of $H^1(\R,\mathbb{C})$, namely $H^1_r$. At the same time, also the stability of the ground state changes if the problem is set in the attractive case or the repulsive one, gaining more stability in the former, as depicted in the following theorems. \\

\begin{thm}[Attractive case, Proposition 2 and Theorem 1 in \cite{foo}]
    \label{thm1}
    Let $\alpha >0$. Then there exists a unique nonnegative minimizer $\phi_{\omega}$ of \eqref{Somega} under the Nehari's constraint. Moreover,
    \begin{itemize}
        \item If $2 < p \leq 6$, then $e^{i\om t}\phi_{\omega}(x)$ is stable in $H^1$ for any $\om \in (\frac{\alpha^2}{4},+\infty)$.
        \item If $p>6$, then there exists $\om_1=\om_1(p,\alpha)>\frac{\alpha^2}{4}$ such that $e^{i\om t}\phi_{\omega}(x)$ is stable in $H^1$ for any $\om \in (\frac{\alpha^2}{4},\om_1)$, and unstable in $H^1$ for any $\om \in (\om_1,+\infty)$. 
    \end{itemize}
\end{thm}
\begin{thm}[Repulsive case, Theorem 1 and 2 in \cite{fj}]
    \label{thm2}
    Let $\alpha <0$. Then there exists a unique nonnegative minimizer $\phi_{\omega}$ of \eqref{Somega} under the Nehari's constraint and among the functions in $H^{1}_{r}$. Moreover,
    \begin{itemize}
        \item If $2 < p \leq 4$, then $e^{i\om t}\phi_{\omega}(x)$ is stable in $H^1_r$ for any $\om \in (\frac{\alpha^2}{4},+\infty)$.
        \item If $4 < p < 6$, then there exists $\om_2=\om_2(p,\alpha)>\frac{\alpha^2}{4}$ such that $e^{i\om t}\phi_{\omega}(x)$ is unstable in $H^1$ for any $\om \in (\frac{\alpha^2}{4},\om_2)$ and stable in $H^1_r$ for any $\om \in (\om_2,+\infty)$.
        \item If $p \geq 6$, then $e^{i\om t}\phi_{\omega}(x)$ is unstable in $H^1$ for any $\om \in (\frac{\alpha^2}{4},+\infty)$. 
    \end{itemize}
\end{thm}

\noindent We note that the value $\om = \frac{\alpha^2}{4}$ corresponds to the frequency of the linear ground state in the attractive case and in particular it represents the threshold after which we can observe the presence of stationary states for the NLSE with delta conditions at the origin. Another remark is that, while in Theorem \ref{thm1} the stability and the instability outcomes hold in $H^1$, in Theorem \ref{thm2} only the instability results are valid in $H^1$ and for the stability ones authors restrict to $H^1_r$.\\\\
These results have been generalized on a star graphs $S_N$ in \cite{acfn2012}, where the search for stationary states has been still conducted both in the attractive and in the repulsive regime. However, the ground state has been identified with the N tail state and characterized as the minimizer of a constrained action functional only in presence of a strong attractive interaction $\alpha^*$. In addition, it has been proved that it is stable in the subcritical and critical regime.\\\\ 
A second family of linear non-Kirchhoff's conditions are the so-called delta prime conditions, introduced in the previous section and defined in \eqref{deltapr}. 
In \cite{an}, we can find a deep investigation about the existence and the orbital stability of ground states using the constrained action functional
\begin{equation*}
    S_{\om, \beta}(u)= \frac{1}{2} \left( ||u'||^2_{L^2(\mathbb{R}_-)}+||u'||^2_{L^2(\mathbb{R}_+)} \right) -\frac{1}{p} ||u||_{L^{p}(\mathbb{R})}^{p} - \frac{1}{2 \beta}|u(0_+)-u(0_-)|^2+\frac{\om}{2} \uLtwor^2.
\end{equation*}
In this work, the authors prove that there exists a critical value $\om^{*}$ for which an interesting bifurcation result occurs. In particular it follows that

\begin{thm}[Theorem 5.3, Proposition 6.11 and Theorem 6.13 in \cite{an}] Let $\beta >0$. Then, there exists $\om^*=\frac{4}{\beta^2}\frac{p}{p-2}$ such that
\begin{itemize}
    \item If $\om \in (\frac{4}{\beta^2},\om^*)$, then there exists a unique ground state, which is odd and orbitally stable for any $p \in (2,6]$.
    \item If $\om \geq \om^*$, then there exist two non-symmetric ground states that are stable if the power nonlinearity does not exceed a critical value $p^*>6$ and become unstable for $p > p^*$. The branch of odd solutions continues to exist at any $\om > \om^*$, but they become a family of orbitally unstable stationary states. 
\end{itemize}
\end{thm}
As it appears from the results, the delta prime conditions give rise to a much richer structure of the family of ground states, including a pitchfork bifurcation with symmetry breaking. In fact, for frequency higher than $\omega^*$, the ground states display no symmetry, making  not possible to reduce the problem to the halfline (contrarily to what happens in the case of a delta, where all ground states are even functions). This higher level of complexity of the whole picture arises from the fact that the energy space is larger, including functions with arbitrary jumps, and no relationship between the positive and negative halflines. Such a connection is restored by the interacting term of the energy.

\noindent More recently, an other type of conditions characterized by a discontinuity has been studied on the real line. They arise from a particular self-adjoint extension of the 1D Laplacian and they are called F\"{u}l\"{o}p-Tsutsui's conditions as in \cite{cheon}. They are defined as
    \begin{equation*}
      \begin{cases}
      \phi(0_+)=\tau \phi(0_-),\\
      \phi'(0_-)- \tau \phi'(0_+)= v \phi(0_-),
\end{cases}    
\end{equation*}
where $\tau \in \R \backslash \{0,1\}$ and $v > 0$. Roughly speaking they can be seen as weighted delta conditions that allow discontinuities at the origin.\\
Some studies have been conducted on these conditions \cite{cheon}, but to the knowledge of the authors, up to now no investigations concerning the existence and the stability of the ground states have been done. To fill the gap and give a more complete review, in the following we present some results to appear \cite{ar} obtained studying the minimization problem for the action functional 
\begin{equation*}
    S_{\om, \tau}(u)= \frac{1}{2} \left( ||u'||^2_{L^2(\mathbb{R}_-)}+||u'||^2_{L^2(\mathbb{R}_+)} \right) -\frac{1}{p} ||u||_{L^{p}(\mathbb{R})}^{p} - \frac{v}{2}|u(0_-)|^2+\frac{\om}{2} \uLtwor^2,
\end{equation*}
on the subset $\Htau : =\{u \in H^1(\mathbb{R}_-) \oplus H^1(\mathbb{R}_+) : u(0_+)=\tau u(0_-)\}$ and under the constraint $I_{\om,\tau}(u)=0$, where $$I_{\om,\tau}(u)= ||u'||^2_{L^2(\mathbb{R}_-)}+||u'||^2_{L^2(\mathbb{R}_+)} - ||u||_{L^{p}(\mathbb{R})}^{p} - v|u(0_-)|^2+ \om \uLtwor^2.$$
The following result proves the existence of ground states for the previous constrained functional
\begin{thm}
    Let $\om >\frac{v^2}{(\tau^2+1)^2}$. Then there exists $u \in \Htau \backslash\{0\}$ that minimizes $S_{\om,\tau}(\cdot)$ and $I_{\om,\tau}(u)=0$. 
\end{thm}
\noindent Even if the proof is quite standard and exploits Banach-Alaoglu's theorem and Brezis-Lieb's lemma in order to obtain a convergence results for the minimizing sequences, a crucial role is played by the following result that allows us to study an equivalent problem.
\begin{prop}
    Let $\om >\frac{v^2}{(\tau^2+1)^2}$. Then
    \begin{align*} 
       d(\om) : & = \inf \{ \, S_{\om,\tau}(u) \, : \, u \in \Htau \backslash \{ 0 \},\, I_{\om,\tau}(u)=0 \} \\
                & = \inf \left\{ \, \frac{p-2}{2p}||u||^p_p \, : \, u \in \Htau \backslash \{ 0 \},\, I_{\om,\tau}(u) \leq 0 \right\}.
\end{align*}
In particular, the two minimization problems are equivalent.
\end{prop}
\noindent Even though these conditions prescribe a discontinuity at the origin, they share the same qualitative behaviour of delta interactions for what concerns the orbital stability of ground states. In particular, since we consider an attractive interaction ($v>0$), we get the following result, analogous to Theorem \ref{thm1} valid for classical delta conditions. More precisely:
\begin{itemize}
    \item If $p \in (2,6]$, then the ground state is stable for any $\om \in \left(\frac{v^2}{(\tau^2+1)^2},+\infty\right)$.
    \item{ [\emph{conjecture}] If $p >6$, then there exists $\bar{\om} > \frac{v^2}{(\tau^2+1)^2}$ such that the ground state is stable for $\om \in \left(\frac{v^2}{(\tau^2+1)^2},\bar{\om}\right)$ and unstable for $\om \in\ (\bar{\om},+\infty)$.}
\end{itemize}      
To conclude, we want to remark that, while the proof of the stability result for $p \in (2,6]$ has been conducted relying on the Grillakis-Shatah-Strauss' theory, the stability conjecture for $p>6$ leans on both numerical simulations and the asymptotical analysis of the behaviour of the $L^2$ norm of the ground state for $\om \in \left( \frac{v^2}{(\tau^2+1)^2}, + \infty \right)$.

\section{Nonlinear delta conditions}
Recently, we started studying the nonlinear Schr\"{o}dinger equation with attractive nonlinear delta interactions at the vertices of the graph. More precisely, starting from the case of the attractive linear delta interaction, it is natural to generalize \eqref{deltaint} replacing the positive number $\alpha$ by $|\phi(v)|^{q-2}$ with $q>2$, getting the condition
\begin{equation}
\label{nldelta}
\sum_{e\succ v} \frac{d\phi_e}{dx_e}(v)=\phi(v)|\phi(v)|^{q-2}\quad \forall v\in \V.
\end{equation}
Such a condition generalizes the model introduced and studied in \cite{at01}, where the effects of nonlinear point interactions are treated, then extended to the three-dimensional setting in \cite{adft03,adft04}, and only recently to space dimension two in \cite{cct,acct1,acct2}. Such models were introduced in order to collect several results coming from theoretical physics and to include concentrated nonlinearities in a new class of mathematically rigorous models. It is worth recalling the application of such models to resonant tunneling \cite{jona,nier}.
An immediate remark is that, differently from all the conditions presented before, \eqref{nldelta} is nonlinear and does not follow from any self-adjoint extension $H$ of the Laplace operator \eqref{laplacian}.\\
As for the linear non-Kirchhoff's conditions, the problem of existence of ground states has been studied only when $\G=\R$ (see \cite{bd}) or when $\G$ is a star graph with $N$ halflines \cite{abd} and not for more general metric graphs yet. In particular, we have looked for global minimizers of the energy functional 
\begin{equation}
\label{Fpq}
F_{p,q}(u,\G)=\frac{1}{2}\int_{\G}|u'|^2\,dx-\frac{1}{p}\int_{\G}|u|^p\,dx-\frac{1}{q}\sum_{v\in \V}|u(v)|^q
\end{equation}
among all the continuous functions $u\in H^{1}(\G)$ satisfying the mass constraint: it is immediate to show that ground states of $\eqref{Fpq}$, if they exist, are solutions of the stationary equation \eqref{stateq} on each edge, are continuous at the vertices by definition and fulfill \eqref{nldelta}. As anticipated in the Introduction, in \eqref{Fpq} we have coexistence of two nonlinearities, the standard and the pointwise one.   
We denote by $\mathcal{F}_{p,q}:[0,+\infty)\to [-\infty,+\infty)$ the function defined as
\begin{equation*}
\label{infFpq}
 \mathcal{F}_{p,q}(\mu):=\inf_{v\in H^1_\mu(S_N)}F_{p,q}(v,S_N).
 \end{equation*}
Due to the presence of two nonlinearities $p$ and $q$ and the validity of Gagliardo-Nirenberg inequalities, if $2<p<6$ and $2<q<4$, then the functional $F_{p,q}$ is bounded from below on $H^{1}_{\mu}(\G)$ for every $\mu>0$ and we are in the so-called subcritical case. When the boundedness from below of $F_{p,q}$ depends on the value of the mass $\mu$, i.e. in the so-called critical cases, it is important to distinguish the single critical case in which $p=6$ and $2<q<4$ or $2<p<6$ and $q=4$ and the doubly critical case in which $p=6$ and $q=4$.\\
\subsection{Subcritical case}
For what concerns the subcritical case, existence results on $\R$ and on $S_{N}$ are very different since sufficient conditions for the compactness of minimizing sequences drastically change.\\
Indeed, one can prove that on $\R$
\begin{equation*}
\inf_{u\in H^{1}_{\mu}(\R)} F_{p,q}(u,\R)<0  \Rightarrow  \text{ a ground state of $F_{p,q}(\cdot,\R)$ exists}
\end{equation*}
and consequently the following theorem holds.
\begin{thm}[Theorem 1.3 in \cite{bd}]
	\label{THM F 1}
	Let $2<p<6$ and $2<q<4$. Then, for every $\mu>0$, there always exists a unique positive ground state of $F_{p,q}(\cdot,\R)$ at mass $\mu$.
\end{thm}
\noindent As one can observe, there is no significant interaction between the two nonlinear terms of the energy since the existence result is the same when the delta interaction is linear.\\
On the other hand, if $\G=S_{N}$, then 
\begin{equation}
\label{scSN}
\inf_{u\in H^{1}_{\mu}(S_{N})} F_{p,q}(S_{N})<\inf_{v\in H^{1}_{\mu}(\R)} E(v,\R)  \Rightarrow  \text{ a ground state of $F_{p,q}(\cdot,S_{N})$ exists}.
\end{equation}
First, notice that, in order to have ground states, the energy $F_{p,q}(\cdot,S_{N})$ has to be smaller than the energy $E(\cdot,\R)$ of the soliton and not only than $0$. More precisely, it can be shown that the existence of a ground state of $F_{p,q}(\cdot, S_{N})$ is equivalent to the existence of a function $u\in H^{1}_{\mu}(S_{N})$ such that $F_{p,q}(u,S_{N})\le E(\phi_{\mu})$. Existence and non-existence results are obtained taking advantage of this equivalence and the interplay between pointwise and standard nonlinearities becomes evident. In particular, if $q<\frac{p}{2}+1$, then existence of ground states holds for small masses and does not hold for large masses and this behaviour is similar to the one described in \cite{acfn2013} in the case of a linear delta interaction at the origin: this suggests that when the nonlinear delta interaction is not too strong, then it has qualitatively the same effect as the linear delta on existence results.  If instead $q>\frac{p}{2}$, then ground states exist for large masses and do not exist for small masses. In both the cases just described, the passage from existence of ground states to non-existence or viceversa identifies a unique threshold value of the mass $\mu^*$ which varies depending on the two powers $p$ and $q$ and on the number of halflines $N$.
When the two nonlinearities are in a perfect balance, that is the case in which $q=\frac{p}{2}+1$, existence and non-existence of ground states depend only on the number of halflines of the star graph and not on the value of the mass.\\
These results are summarized in the following theorem, which is in preparation \cite{abd}.
 \begin{thm}
Let $2<p<6$, $2<q<4$ and $N$ the number of halflines of the star graph.\\
If $q<\frac{p}{2}+1$, then there exists $\mu^*(p,q,N)>0$ such that
\begin{itemize}
\item{if $\mu\le\mu^*(p,q,N)$, then there exists a ground state of \eqref{Fpq} at mass $\mu$,}
\item{if $\mu>\mu^*(p,q,N)$, then $\F_{p,q}(\mu)$ is not attained.}
\end{itemize}
On the contrary, if $q>\frac{p}{2}+1$, then there exists $\mu^*(p,q,N)>0$ such that
\begin{itemize}
\item{if $\mu<\mu^*(p,q,N)$, then $\F_{p,q}(\mu)$ is not attained,}
\item{if $\mu\ge\mu^*(p,q,N)$, then there exists a ground state of \eqref{Fpq} at mass $\mu$.}
\end{itemize}
If instead $q=\frac{p}{2}+1$, then there exists $N^*(p)\ge 3$ such that 
   \begin{itemize}
   \item{if $N\le N^*(p)$, then for every $\mu>0$ there exists $u\in H^1_\mu(S_N)$ such that $F_{p,q}(u)=\F_{p,q}(\mu)$,}
   \item{if $N>N^*(p)$, then for every $\mu>0$ no ground state of \eqref{Fpq} at mass $\mu$ exists.}
   \end{itemize}
   \end{thm}
   
   \subsection{Critical cases}
Existence of ground states in critical cases has been studied only when $\G=\R$. The first theorem deals with the cases in which only one power is critical. 
  \begin{thm}[Theorem 1.4 in \cite{bd}]
	\label{THM F 2}
	Let $\mu>0$. 
	\begin{itemize}
		\item[(i)] If $p=6$ and $2<q<4$, then there exists a unique positive ground state at mass $\mu$ if and only if $\mu<\frac{\sqrt{3}}{2} \pi$, and
		\begin{equation}
		\label{eq_THM F 2}
		\begin{cases}
		-\infty<\F_{6,q}(\mu)<0 & \text{if }\mu<\frac{\sqrt{3}}2 \pi\\
		\F_{6,q}(\mu)=-\infty & \text{if }\mu\geq\frac{\sqrt{3}}2 \pi\,.
		\end{cases}
		\end{equation}
		\item[(ii)] If $2<p<6$ and $q=4$, then there exists a unique positive ground state at mass $\mu$ if and only if $\mu<2$
		\begin{equation}
		\label{eq_THM F 3}
		\begin{cases}
		-\infty<\F_{p,4}(\mu)<0 & \text{if }\mu<2\\
		\F_{p,4}(\mu)=-\infty & \text{if }\mu\geq2\,.
		\end{cases}
		\end{equation}
	\end{itemize}
\end{thm}
\noindent These particular regimes show the interplay between a subcritical and a critical power nonlinearity: indeed, while the ground state level moves on from $0$ to $-\infty$ in correspondence of a threshold value of the mass as usual in critical cases, the presence of the subcritical power ensures existence of ground states for all the masses under the critical mass, highlighting an important difference with what one expects in critical cases. \\

\noindent The last result concerns the doubly critical case, where simultaneously $p=6$ and $q=4$. Here we recover the typical structure of a purely critical setting, with the ground state energy level lifting from 0 to $-\infty$ when exceeding a critical value of the mass and solutions existing only at the threshold. A quite remarkable feature due to the interaction between the two nonlinearities is given by the fact that the critical mass \eqref{mustar} is lower than the critical masses $\frac{\sqrt{3}}2\pi$ and $2$ for the standard and pointwise nonlinearity.
\begin{thm}[Theorem 1.5 in \cite{bd}]
	\label{THM F 4}
	The functional $F_{6,4}(\cdot,\R)$ admits ground states only at mass 
	\begin{equation}
	\label{mustar}
	\mu^*:=\sqrt{3}\left(\frac{\pi}{2}-\arcsin\left(\sqrt{\frac{3}{7}}\right)\right)
	\end{equation}
	and
	\begin{equation}
	\label{eq_THM F 4}
	\F_{6,4}(\mu)=\begin{cases}
	0 & \text{if }\mu\leq\mu^*\\
	-\infty & \text{if }\mu>\mu^*\,.
	\end{cases}
	\end{equation}
\end{thm}

\section{Conclusion and perspectives}
We reviewed on recent developments and announced new results on the the existence of Ground States for the Nonlinear Schr\"odinger Equation on infinite star graphs, stressing possible applications to the description of the behaviour of Bose-Einstein condensates. From the physical point of view we are still at a preliminary step, in which the correct conditions are to be singled out, but there are hints that the correct conditions could not be free, or Kirchhoff's. On the other hand, such a topic has already shown a very rich potential from a mathematical point of view, underlying how the interplay between a point interaction and a standard power nonlinearity can give rise to phenomena that are new, in the sense that they do not occur neither in the standard Nonlinear Schr\"odinger equation, nor in the linear one with point interactions. \\
\noindent
The most spectacular among them is the occurrence of bifurcations, for which a non-trivial energy space is needed, a requirement that is fulfilled by specific (in particular, discontinuous) boundary conditions at the vertex, e.g. delta prime or F\"ul\"op-Tsutsui.
Further phenomena to be still completely investigated arise when the vertex interaction is nonlinear too, so that there is competition between the two nonlinearities and possibly coexistence of two criticalities. In our results it is shown how this competition is ruled, and in which cases one nonlinearity prevails on the others.\\
\noindent
As it is natural, a prosecution of the research will be concerned with the deepening of the resarch in this direction. Another one, completely to be explored, will deal with considering more-dimensional analogues of the graphs, the so-called hybrids, with components of different dimensions to be glued together. Preliminary discussions and formal computations made on this, seem to show that the phenomenology given by the nonlinear hybrids could be still richer than the one exhibited by nonlinear graphs.


\begin{thebibliography}{99}

\bibitem{abd} Adami R., Boni F., Dovetta S., {\em Ground states of a doubly nonlinear Schrödinger equation on star graphs}, in preparation.



\bibitem{acfn2012} Adami R., Cacciapuoti C., Finco D., Noja D., {\em Stationary states of NLS on star graphs}, EPL (Europhysics Letters) 100 (1), 10003 (2012).

\bibitem{acfn2013} Adami R., Cacciapuoti C., Finco D., Noja D., {\em Constrained energy minimization and orbital stability for the NLS equation on a star graph}, Ann. I. H. Poincaré – AN (2013).

\bibitem{acfn2016} Adami R., Cacciapuoti C., Finco D., Noja D., {\em Stable standing waves for a NLS on star graphs as local minimizers of the constrained energy}, Journal of Differential Equations 260 (10), 7397--7415 (2016).

\bibitem{acct1} Adami R., Carlone R., Correggi M., Tentarelli L., {\em Blow-up for the pointwise NLS in dimension two: Absence of critical power}, J. Diff. Eq. {\bf 269} (1), 1--37 (2020).

\bibitem{acct2} Adami R., Carlone R., Correggi M., Tentarelli L.,
{\em Stability of the standing waves of the concentrated NLS in dimension two}, preprint arXiv:2001.03969, (2020).
		


\bibitem{adft03}	Adami R., Dell’Antonio G., Figari R., Teta A.: {\em The Cauchy problem for the Schrödinger equation in dimension three with concentrated nonlinearity}. Ann. I. H. Poincaré. {\bf20}, 477–500 (2003).

\bibitem{adft04} Adami R., Dell’Antonio G., Figari R., Teta A., {\em Blow-up solutions for the Schrödinger equation in dimension three with a concentrated nonlinearity}. Ann. I. H. Poincaré. {\bf 21}, 121–137 (2004).

\bibitem{ad} Adami R., Dovetta S., {\em One-dimensional versions of three-dimensional system: Ground states for the NLS on the spatial grid}, Rend. Mat. {\bf 39} (2018), 181--194.

\bibitem{adst}
		Adami R., Dovetta S., Serra E., Tilli P., 
		{\em Dimensional crossover with a continuum of critical exponents for NLS on doubly periodic metric graphs}, Analysis \& PDE, Vol. 12 (2019), No. 6, 1597–-1612.
		
\bibitem{adr} Adami R., Dovetta S., Ruighi A., {\em Quantum graphs and dimensional crossover: the honeycomb}, COMMUNICATIONS IN APPLIED AND INDUSTRIAL MATHEMATICS. - ISSN 2038-0909. - 10:1, 109--122 (2019).		

\bibitem{an} Adami R., Noja D., {\em Stability and symmetry breaking bifurcation for the ground states of a NLS equation with a $\delta'$ interaction}, Communications in Mathematical Physics 318 (1), 247--289 (2013).

\bibitem{anv} Adami R., Noja D., Visciglia N., {\em Constrained energy minimization and ground states for NLS with point defects}, Discrete Contin. Dyn. Syst. Ser. B., 18 (2013), 1155–1188.

\bibitem{ar} Adami R., Ruighi A., {\em Discontinuous ground states for the NLSE on $\R$ with a F\"ul\"op-Tsutsui $\delta$ type interaction}, in preparation.

\bibitem{ast17}
Adami R., Serra E., Tilli P. 
{\em Negative energy ground states for the $L^2$--critical NLSE on metric graphs}, 
Comm. Math. Phys. {\bf 352} (2017), no. 1, 387--406.
		
\bibitem{ast15}
Adami R., Serra E., Tilli P.:
{\em NLS ground states on graphs}.
Calc. Var. and PDEs {\bf 54} (2015) no. 1, 743--761.

\bibitem{ast-jfa}
		Adami R., Serra E., Tilli P., 
		{\em Threshold phenomena and existence  for NLS ground states on graphs}, J. Funct. An. \textbf{271(1)}, 201--223 (2016).

\bibitem{at01} Adami R., Teta A., {\em A class of nonlinear Schr\``odinger equation with concentrated nonlinearity}, J. Funct. An. {\bf 180},  148--175 (2001).

\bibitem{alb} Albeverio S., Gesztesy F., Hoegh-Krohn R., Holden H., {\em Solvable Models in Quantum
Mechanics}, Springer, New York, (1988).

\bibitem{alimehmeti} Ali Mehmeti, F. Nonlinear waves in Networks, Wiley VCH, 1994.

\bibitem{kuchment} Berkolaiko G., Kuchment P. Introduction to Quantum Graphs, Mathematical Surveys and Monographs, AMS, 2012.



\bibitem{bd} Boni F., Dovetta S., {\em Ground states for a doubly nonlinear Schrödinger equation in dimension one}, arXiv:1907.07926 [math.AP] (2019).

\bibitem{bose} Bose S.N., {\em Plancks Gesetz und Lichtquantenhypothese}, Zeit. f\"ur Physik,  
{\bf 26} (1924), 178--181.

\bibitem{ccdt19} Carlone R., Correggi M., Finco D., Teta A., {\em A Quantum Model of Feshbach Resonances}, Ann. H. Poincar\'e {\bf 20} (2019), 2899--2935.  



\bibitem{cct} Carlone R., Correggi M., Tentarelli L., {\em Well-posedness of the two-dimensional nonlinear Scr\"odinger equation with concentrated nonlinearity}, Ann. Inst. H. Poincar\'e An. Non Lin {\bf 36}, 1, 257--294 (2019).
 
 


\bibitem{lions}
 Cazenave T., Lions P.-L.
 {\em Orbital stability of standing waves for some nonlinear Schr\"odinger equations}. Commun. Math. Phys. {\bf 85} (1982), no. 4, 549--561. 
 
\bibitem{c15} Caudrelier V., {\em On the inverse scattering method for integrable PDEs on a star graph}, Commun. Math. Phys. {\bf 338} (2) (2015), 893--917.

\bibitem{cmr1} Caudrelier V., Mintchev M., Ragoucy E., {\em Solving the quantum non-linear Schr\"odinger
equation with delta-type impurity}, J. Math. Phys. 46 (2005) 042703.

\bibitem{cmr2} Caudrelier V., Mintchev M., Ragoucy E., {\em The quantum non-linear Schr\"odinger
model with point-like defect}, J.Phys. A37 (2004) L367-L376.

\bibitem{cr} Caudrelier V., Ragoucy E., {\em Direct computation of scattering matrices for general
quantum graphs}, Nucl. Phys. B828 (2010) 515-535.

\bibitem{cheon} Cheon T., Turek O., {\em Fulop–Tsutsui interactions on quantum graphs},Physics Letters A {\bf 374} (41) (2010), 4212--4221.
 
\bibitem{dovetta-per} Dovetta S., {\em Mass-constrained ground states of the stationary NLSE on periodic metric graphs}, Nonlinear Differ. Equ. Appl. {\bf 26} (2019), n. 30.


\bibitem{dst}
Dovetta S., Serra E., Tilli P., {\em NLS  ground  states  on  metric  trees:existence results and open questions}, arXiv preprint arXiv:1905.00655, (2019).


\bibitem{dt}
Dovetta S., Tentarelli L.,
{\em $L^2$--critical NLS on noncompact metric graphs with localized nonlinearity: topological and metric features}, Calc Var. PDE {\bf 58} (3) (2019), n. 108.


\bibitem {einstein} Einstein A., {\em Quantentheorie des einatomigen idealen Gases}, Sitz.  Preus. Akad. Wiss., {\bf 1} (1925), 3.

\bibitem{fj} Fukuizumi R., Jeanjean L., {\em Stability of standing waves for a nonlinear Schr\"{o}dinger equation with a repulsive Dirac delta potential}, Disc. Cont. Dyn. Syst. (A), {\bf 21}, 129--144 (2008).

\bibitem{foo} Fukuizumi R., Otha M., Ozawa T., {\em Nonlinear Schr\"{o}dinger equation with a point defect}, Ann. IHP, Analyse non linéaire, {\bf 25}, 837--845 (2008). 


\bibitem{ft} F\"ul\"op T., Tsutsui I., {\em A free particle on a circle with point interaction}, Phys. Lett. A 264 (5), 366--374 (2000).


\bibitem{ghw} Goodman R.H., Holmes P.J., Weinstein M.I., {\em Strong NLS soliton-defect interaction}, Physica D, {\bf 192}, 215--248 (2004).

\bibitem{gss1} Grillakis M., Shatah J., Strauss W., {\em Stability theory of solitary waves in the presence of symmetry - I}, J. Func. An., {\bf 74}, 160--197 (1987).

\bibitem{gss2} Grillakis M., Shatah J., Strauss W., {\em Stability theory of solitary waves in the presence of symmetry - II}, J. Func. An., {\bf 94}, 308-348 (1990).

\bibitem{hmz} Holmer J., Marzuola J., Zworski M., {\em Fast soliton scattering by delta impurities}, Comm. Math. Phys., {\bf 274}, 187--216 (2007).


\bibitem{jona} Jona Lasinio G., Presilla C., Sj\"ostrand P., {\em On Schr\"odinger Equations with Concentrated Nonlinearities}, Ann. of Phys. {\bf 240} (1995), 1--21. 

\bibitem{kp2} Kairzhan A., Pelinovsky D.E., {\em Nonlinear instability of half-solitons on star-graphs}, J. Diff. Eq. {\bf 264} (2018), 7357--7383.


\bibitem{kp} Kairzhan A., Pelinovsky D.E., {\em Spectral stability of shifted states on star graphs}, J. Phys. A Math. Theor. {\bf 51} (9) (2018), 095203.


\bibitem{kpg} Kairzhan A., Pelinovsky D.E., Goodman R. H., {\em Drift of spectrally stable shifted states on star graphs}, SIAM J. Appl. Dyn. Sys. {\bf 18} (4) (2019), 1723--1755.

\bibitem{kost} Kostrykin V., Schrader R., {\em Kirchhoff's Rule for Quantum Wires}, J. Phys. A32 (1999) 595-630.




\bibitem{coz} Le Coz S., Fukuizumi R., Fibich G., Ksherim B., Sivan Y., {\em Instability of bound states of a nonlinear Schr\"{o}dinger equation with a Dirac potential}, Physica D {\bf 237}, No. 8, 1103-1128 (2008).

\bibitem{nier} Nier, F., {\em The dynamics of some open quantum systems with short-range nonlinearities}, Nonlinearity {\bf 11}, 4, 1127--1172 (1998).

\bibitem{noja14} Noja D., {\em Nonlinear Schr\"odinger on graphs: recent results and open problems}, Phil. Trans. R. Soc. A, {\bf 372} (2014), 20130002.

\bibitem{serratentarellly}
Serra E., Tentarelli L.
{\em Bound states of the NLS equation on metric graphs with localized nonlinearities}.
J. Diff. Eq. {\bf 260} (2016), no. 7, 5627--5644.


\bibitem{sbmnu} Sobirov Z. A., Babajanov D., Matrasulov D., Nakamura K., Uecker H., {\em Sine-Gordon soliton in networks: scattering and transmission at vertices}, Europhys Lett. {\bf 115} (2016), 50002. 


\bibitem{ssms} Sobirov Z. A., Matrasulov D., Sawada S., Nakamura K., {\em Integrable nonlinear Schr\"odinger equation on simple networks: connection formula at vertices}, Phys. Rev. E {\bf 81} (6-2) (2010), 066602.








\bibitem{t}
Tentarelli L.,
{\em NLS ground states on metric graphs with localized nonlinearities},
\emph{J. Math. Anal. Appl.} {\bf 433} (2016), no. 1, 291--304.

\bibitem{ftc} Tsutsui I., F\"ul\"op T,  Cheon T., {\em Connection conditions and the spectral family under singular potentials}, J. Phys. A Math. Gen. {\bf 36} 1, 275 (2002).




\bibitem{zakharov}
Zakharov V.E., Shabat B., {\em Exact Theory of Two--Dimensional Self--Focusing and One--Dimensional Self--Modulation of Waves in Nonlinear Media},Soviet Phys. JETP {\bf34} (1) (1972), 62-–71.

\end{thebibliography}
\end{document}